\documentclass[11pt,leqno]{article}
\usepackage{amssymb}
\usepackage{amsfonts}
\newtheorem{prop}{Proposition}[section]
\newtheorem{teo}[prop]{Theorem}
\newtheorem{remark}[prop]{Remark}

\newcommand{\trace}{\textnormal{trace}}

\date{}
\begin{document}
\author{D. FETCU and C. ONICIUC\thanks{The authors were partially supported by the Grant CEEX, ET, 5871/2006, Romania.}
\\ \\ \textit{Dedicated to Professor Vasile Oproiu on his 65-th
birthday}}
\title{Explicit formulas for biharmonic submanifolds in non-Euclidean 3-spheres}
\date{}
\maketitle
\begin{abstract}
We obtain the parametric equations of all biharmonic Legendre
curves and Hopf cylinders in the 3-dimensional unit sphere endowed
with the modified Sasakian structure defined by Tanno.
\end{abstract}

\section{Introduction}
\setcounter{equation}{0}

\textit{Biharmonic maps} $\phi:(M,g)\rightarrow(N,h)$ between
Riemannian manifolds are the critical points of the \textit{bienergy
functional} $E_{2}(\phi)=\frac{1}{2}\int_{M}|\tau(\phi)|^{2}\
\nu_{g}$ and represent a natural generalization of the well-known
\textit{harmonic maps}, the critical points of the \textit{energy
functional} $E(\phi)=\frac{1}{2}\int_{M}|d\phi|^{2}\ \nu_{g}$. The
bienergy functional has been suggested since 1964 by Eells and
Sampson in their famous paper \cite{Eells}. The bienergy functional
can be viewed as the analogue for maps of the Willmore functional
for Riemannian immersions. In a different setting, Chen defined the
biharmonic submanifolds in an Euclidean space. If we apply the
characterization formula of biharmonic maps to Riemannian immersions
into Euclidean spaces, we recover Chen's notion of biharmonic
submanifold.

The Euler-Lagrange equation for the energy functional is
$\tau(\phi)=0$, where $\tau(\phi)=\trace\nabla d\phi$ is the
tension field, and the Euler-Lagrange equation for the bienergy
functional was derived by Jiang in \cite{Jiang}:
$$
\begin{array}{cl}
\tau_{2}(\phi)&=-\Delta\tau(\phi)-\trace\
R^{N}(d\phi,\tau(\phi))\tau(\phi)=\\ \\
&=0.\end{array}
$$
Since any harmonic map is biharmonic, we are interested in
non-harmonic biharmonic maps, which are called
\textit{proper-biharmonic}.

There are several classification results of proper-biharmonic
submanifolds in 3-dimensional spaces. The proper-biharmonic
submanifolds in $\mathbb{R}^{3}$, $\mathbb{S}^{3}$ and $N^{3}(-1)$,
the space of constant negative sectional curvature $-1$, were
completely classified in \cite{Chen1}, \cite{Cad} and \cite{Cad1},
respectively. Then, the next step was to classify proper-biharmonic
submanifolds in a 3-dimensional space of non-constant sectional
curvature. In \cite{COP}, \cite{Cad2} and \cite{CIL} the authors
classified the proper-biharmonic curves in the Heisenberg group and
Cartan-Vranceanu 3-dimensional spaces, while, in \cite{Ino},
Inoguchi classified the proper-biharmonic Legendre curves and Hopf
cylinders in a 3-dimensional Sasakian space form $M^{3}(c)$.

For a general account of biharmonic maps see
\cite{MontaldoOniciuc} and \textit{The bibliography of biharmonic
maps} \cite{bibl}.

The goal of this paper is to obtain explicitly the parametric
equations of proper-biharmonic Legendre curves and Hopf cylinders
in $\mathbb{S}^{3}$ with the modified Sasakian structure defined
by Tanno, by using the results in \cite{Ino}. The techniques used
in our paper are introduced in \cite{BB} and \cite{BBK}.

\noindent \textbf{Conventions.} We work in the $C^{\infty}$
category, that means manifolds, metrics, connections and maps are
smooth. The Lie algebra of the vector fields on $M$ is denoted by
$C(TM)$.

\section{Preliminaries}
\setcounter{equation}{0}

\subsection{Contact manifolds}

A \textit{contact metric structure} on a manifold $M^{2n+1}$ is
given by $(\varphi,\xi,\eta,g)$, where $\varphi$ is a tensor field
of type $(1,1)$ on $M$, $\xi$ is a vector field on $M$, $\eta$ is
an 1-form on $M$ and $g$ is a Riemannian metric, such that
$$
\left\{\begin{array}{lll} \varphi^{2}=-I+\eta\otimes\xi,\ \
\eta(\xi)=1,\\ \\ g(\varphi X,\varphi Y)=g(X,Y)-\eta(X)\eta(Y),\ \
\forall X,Y\in C(TM)\\ \\ g(X,\varphi Y)=d\eta(X,Y),\ \ \forall
X,Y\in C(TM)\end{array}\right..
$$
A contact metric structure $(\varphi,\xi,\eta,g)$ is
\textit{Sasakian} if it is normal, that is
$$
N_{\varphi}+2d\eta\otimes\xi=0,
$$
where
$$
N_{\varphi}(X,Y)=[\varphi X,\varphi Y]-\varphi \lbrack \varphi
X,Y]-\varphi \lbrack X,\varphi Y]+\varphi ^{2}[X,Y],\ \ \forall
X,Y\in C(TM),
$$
is the Nijenhuis tensor field of $\varphi$, or, equivalently, if
$$
(\nabla_{X}\varphi)(Y)=g(X,Y)\xi-\eta(Y)X,\ \ \forall X,Y\in
C(TM).
$$

Let us consider a Sasakian manifold $(M,\varphi,\xi,\eta,g)$. The
sectional curvature of a 2-plane generated by $X$ and $\varphi X$,
where $X$ is an unit vector orthogonal to $\xi$, is called
\textit{$\varphi$-sectional curvature} determined by $X$. A Sasakian
manifold with constant $\varphi$-sectional curvature $c$ is called a
\textit{Sasakian space form} and it is denoted by $M(c)$.

\noindent The curvature tensor field of a Sasakian space form
$M(c)$ is given by
$$
\begin{array}{ll}
R(X,Y)Z=&\frac{c+3}{4}\{g(Z,Y)X-g(Z,X)Y\}+\frac{c-1}{4}\{\eta(Z)\eta(X)Y-\\ \\
&-\eta(Z)\eta(Y)X+g(Z,X)\eta(Y)\xi-g(Z,Y)\eta(X)\xi+\\
\\&+g(Z,\varphi Y)\varphi X-g(Z,\varphi X)\varphi
Y+2g(X,\varphi Y)\varphi Z\}.
\end{array}
$$

A contact metric manifold $(M,\varphi,\xi,\eta,g)$ is called
\textit{regular} if for any point $p\in M$ there exists a cubic
neighborhood of $p$ such that any integral curve of $\xi$ passes
through the neighborhood at most once, and \textit{strictly regular}
if all integral curves are homeomorphic to each other.

\noindent Let $(M,\varphi,\xi,\eta,g)$ be a regular contact metric
manifold. Then the orbit space $\bar{M}$ can be naturally organized
as a manifold and, moreover, if $M$ is compact then $M$ is a
principal circle bundle over $\bar{M}$ (the Boothby-Wang Theorem).
In this case the fibration $\pi:M\rightarrow\bar{M}$ is called
\textit{the Boothby-Wang fibration}. A very known example of a
Boothby-Wang fibration is the Hopf fibration
$\pi:\mathbb{S}^{2n+1}\rightarrow \mathbb{C}P^{n}$.

We recall the following result obtained by Ogiue
\begin{teo}[\cite{Ogiue}] Let $(M,\varphi,\xi,\eta,g)$ be a strictly
regular Sasakian manifold. Then $\bar{M}$ can be organized as a
K\"{a}hler manifold. Moreover, if $(M,\varphi,\xi,$ $\eta,g)$ is a
Sasakian space form $M(c)$, then $\bar{M}$ has constant sectional
holomorphic curvature $c+3$.
\end{teo}

\subsection{General results}

In \cite{Cad} all biharmonic curves and surfaces in the
3-dimensional Euclidean unit sphere $(\mathbb{S}^{3},g_{0})$,
where $g_{0}$ is the usual metric, were classified:

\begin{teo}[\cite{Cad}] Let $M^{m}$ be a proper-biharmonic
submanifold of $\mathbb{S}^{3}$. We have

a) if $m=1$, then $M$ is either the circle of radius
$\frac{1}{\sqrt{2}}$ or a geodesic of the torus
$\mathbb{S}^{1}\Big(\frac{1}{\sqrt{2}}\Big)\times
\mathbb{S}^{1}\Big(\frac{1}{\sqrt{2}}\Big)\subset \mathbb{S}^{3}$
with the slope different from $\pm 1$;

b) if $m=2$, then $M$ is the hypersphere
$\mathbb{S}^{2}\Big(\frac{1}{\sqrt{2}}\Big)\subset
\mathbb{S}^{3}$.
\end{teo}

We can think $(\mathbb{S}^{3},g_{0})$ as a Sasakian space form with
constant $\varphi_{0}$-sectional curvature 1. We know that in the
geometry of Sasakian manifolds an important role is played by the
integral submanifolds and Hopf cylinders. A natural question is
whether biharmonic Legendre curves and Hopf cylinders in
$(\mathbb{S}^{3},\varphi_{0},\xi_{0},\eta_{0},g_{0})$ exist. Since
the torsion of a Legendre curve in
$(\mathbb{S}^{3},\varphi_{0},\xi_{0},\eta_{0},g_{0})$ is 1 and a
Hopf cylinder is flat, it follows that in
$(\mathbb{S}^{3},\varphi_{0},\xi_{0},$ $\eta_{0},g_{0})$ do not
exist proper-biharmonic Legendre curves and proper-biharmonic Hopf
cylinders.

The next step is the study of the existence of biharmonic Legendre
curves and of biharmonic Hopf cylinders in Sasakian space forms with
constant $\varphi$-sectional curvature $c\neq 1$. Inoguchi gave the
following classification:

\begin{teo}[\cite{Ino}]\label{t1.1}
Let $M^{3}(c)$ be a Sasakian space form of constant
$\varphi$-sectional curvature $c$ and let $\gamma:I\rightarrow M$
be a Legendre curve parametrized by arc length. We have

a) if $c\leqslant 1$, then $\gamma$ is biharmonic if and only if it
is geodesic;

b) if $c>1$, then $\gamma$ is proper-biharmonic if and only if it
is a helix with the curvature $\bar\kappa^{2}=c-1$.
\end{teo}

\begin{teo}[\cite{Ino}]\label{t1.2}
Let $M^{3}(c)$ be a Sasakian space form of constant
$\varphi$-sectional curvature $c$ and $S_{\bar\gamma}$ a Hopf
cylinder, where $\bar\gamma$ is a curve in the orbit space of
$M^{3}(c)$, parametrized by arc length. We have

a) if $c\leqslant 1$, then $S_{\bar\gamma}$ is biharmonic if and
only if it is minimal;

b) if $c>1$, then $S_{\bar\gamma}$ is proper-biharmonic if and
only if the curvature $\bar\kappa$ of $\bar\gamma$ is constant
$\bar\kappa^{2}=c-1$.
\end{teo}

Let $M^{3}(c)$ be a Sasakian space form with constant
$\varphi$-sectional curvature $c>1$, and let
$\bar\gamma:I\rightarrow \bar{M}$ be a curve parametrized by arc
length. We denote $f_{1}=\bar\gamma '$ and $f_{2}=\bar{J}f_{1}$ the
Frenet frame field along $\bar\gamma$, where $\bar{J}$ is the almost
complex structure on $\bar{M}$. Consider
$S_{\bar\gamma}=\pi^{-1}(\bar\gamma)$ the Hopf cylinder
corresponding to $\bar\gamma$ and assume that it is biharmonic, that
is $\bar{k}^{2}=c-1$. Since $[\xi,f_{1}^{H}]=0$, $g(\xi,\xi)=1$,
$g(f_{1}^{H},f_{1}^{H})=1$, we can choose a local chart $x=x(u,v)$
such that $\xi=x_{v}$ and $f_{1}^{H}=x_{u}$. The parametric curves
$u\rightarrow x(u,v_{0})$ are proper-biharmonic Legendre curves
parametrized by arc length and geodesics in $S_{\bar\gamma}$.
Moreover, $u\rightarrow x(u,v_{0})$ are the only geodesics of
$S_{\bar\gamma}$ which are proper-biharmonic in $M^{3}(c)$.

\begin{prop}
Let $M^{3}(c)$ be a Sasakian space form with constant
$\varphi$-sectional curvature $c>1$ and $S_{\bar\gamma}$ a
biharmonic Hopf cylinder. Let $\gamma:I\rightarrow S_{\bar\gamma}$
be a geodesic parametrized by arc length and
$\textnormal{\textbf{i}}:S_{\bar{\gamma}}\rightarrow M^{3}(c)$ the
canonical inclusion. Then
$\textnormal{\textbf{i}}\circ\gamma:I\rightarrow M^{3}(c)$ is
proper-biharmonic if and only if $\gamma'=\pm f_{1}^{H}$.
\end{prop}
\begin{pf} Let $\gamma:I\rightarrow S_{\bar\gamma}$ be a geodesic
parametrized by arc length. Then $\gamma'=c_{1}\xi+c_{2}f_{1}^{H}$,
where $c_{1}$ and $c_{2}$ are real constants such that
$c_{1}^{2}+c_{2}^{2}=1$ and $f_{1}^{H}$ is the horizontal lift of
$f_{1}$. After a straightforward computation one obtains
$$
\tau(\textbf{i}\circ\gamma)=c_{2}(c_{2}\bar\kappa-2c_{1})f_{2}^{H}\
\ \hbox{and}\ \
\tau_{2}(\textbf{i}\circ\gamma)=2c_{1}c_{2}^{2}(c_{2}\bar\kappa-2c_{1})\bar\kappa
f_{2}^{H},
$$
where $\bar\kappa^{2}=c-1$. Hence $\tau(\textbf{i}\circ\gamma)\neq
0$ and $\tau_{2}(\textbf{i}\circ\gamma)=0$ if and only if $c_{1}=0$,
that is $\gamma'=\pm f_{1}^{H}$.

\hskip11cm $\Box$

\end{pf}

It can be easily proved that, if $\gamma:I\rightarrow M^{3}(c)$ is
a biharmonic Legendre curve, then it is a geodesic in the
biharmonic Hopf cylinder $x(s,t)=\phi_{\gamma(s)}(t)$, where
$\{\phi_{t}\}$ is the uniparametric group of $\xi$. From now on we
will use the parameters $(u,v)$ instead of $(s,t)$.

In the following we will choose the 3-sphere $\mathbb{S}^{3}$ with
modified Sasakian structure as a model for the complete, simply
connected Sasakian space form with constant $\varphi$-sectional
curvature $c>1$, and we will find the explicit equations of
biharmonic Legendre curves and of biharmonic Hopf cylinders, viewed
as submanifolds of $\mathbb{R}^{4}$. These equations are given by
Theorem \ref{t2} and Theorem \ref{t3}, respectively.

\section{Proper-biharmonic Legendre curves} \setcounter{equation}{0}

Let $\mathbb{S}^{3}=\{z\in\mathbb{C}^{2}: \|z\|=1\}$ be the unit
3-dimensional sphere endowed with its standard metric field $g_{0}$.
Consider the following structure tensor fields on $\mathbb{S}^{3}$:
$\xi_{0}=-Jz$ for each $z\in \mathbb{S}^{3}$, where $J$ is the usual
almost complex structure on $\mathbb{C}^{2}$ defined by
$Jz=(-y_{1},-y_{2},x_{1},x_{2})$ for $z=(x_{1},x_{2},y_{1},y_{2})$,
and $\varphi_{0}=s\circ J$, where $s:T_{z}\mathbb{C}^{2}\rightarrow
T_{z}\mathbb{S}^{3}$ denotes the orthogonal projection. Endowed with
these tensors, $\mathbb{S}^{3}$ becomes a Sasakian space form with
$\varphi_{0}$-sectional curvature 1 (see \cite{Blair}). Now we
consider the deformed structure
$$
\eta=a\eta_{0}, \ \ \xi=\frac{1}{a}\xi_{0},\ \
\varphi=\varphi_{0},\ \ g=a g_{0}+a(a-1)\eta_{0}\otimes\eta_{0},
$$
where $a$ is a positive constant. Such a deformation is called a
$\mathcal{D}$-homothetic deformation, since the two metrics
restricted to the contact subbundle $\mathcal{D}$ are homothetic,
and it was introduced in \cite{Tanno}. The deformed structure
$(\varphi,\xi,\eta,g)$ is still a Sasakian structure and
$(\mathbb{S}^{3},\varphi,\xi,\eta,g)$ is a Sasakian space form
with constant $\varphi$-sectional curvature $c=\frac{4}{a}-3$ (see
\cite{BBK}). From now on we assume that $a\in (0,1)$, that is
$c>1$.

\begin{teo}\label{t2}
Let $\gamma:I\rightarrow (\mathbb{S}^{3},\varphi,\xi,\eta,g)$ be a
Legendre curve parametrized by arc length. Then it is
proper-biharmonic if and only if, as a curve in $\mathbb{R}^{4}$
\begin{equation}\label{ec1}
\begin{array}{cl}
\gamma(s)=&\sqrt{\frac{B}{A+B}}\cos(As)e_{1}-\sqrt{\frac{B}{A+B}}\sin(As)Je_{1}+\\
\\&+\sqrt{\frac{A}{A+B}}\cos(Bs)e_{3}+\sqrt{\frac{A}{A+B}}\sin(Bs)Je_{3},\end{array}
\end{equation}
where $\{e_{1},e_{3}\}$ is an orthonormal system of constant
vectors in the Euclidian space $(\mathbb{R}^{4},\langle,\rangle)$,
with $e_{3}$ orthogonal to $Je_{1}$, and
\begin{equation}\label{2}
\left\{\begin{array}{ll}A=\sqrt{\frac{3-2a-2\sqrt{(a-1)(a-2)}}{a}}\\
\\ B=\sqrt{\frac{3-2a+2\sqrt{(a-1)(a-2)}}{a}}\end{array}\right..
\end{equation}
\end{teo}

\begin{pf} Let us denote by $\nabla$, $\dot{\nabla}$ and
$\widetilde{\nabla}$ the Levi-Civita connections on
$(\mathbb{S}^{3},g)$, $(\mathbb{S}^{3},g_{0})$ and
$(\mathbb{R}^{4},\langle,\rangle)$, respectively.

Let $\gamma:I\rightarrow(\mathbb{S}^{3},g)$ be a proper-biharmonic
Legendre curve parametrized by arc length and let $T=\gamma'$ be
the unit tangent vector field along the curve. Using Theorem
\ref{t1.1} we have
$$
\bar{g}(T,T)=1,\ \ g(T,\xi)=0,\ \ \nabla_{T}T=\kappa\varphi T,\ \
\hbox{and}\ \ \kappa^{2}=c-1.
$$

In order to find the explicit parametrization of $\gamma$ as  a
curve in $\mathbb{R}^{4}$ we shall characterize $\gamma$ as the
solution of a certain fourth order ODE with constant coefficients.
For this purpose the expressions of $\widetilde{\nabla}_{T}T$,
$\widetilde{\nabla}_{T}\varphi T$ and $\widetilde{\nabla}_{T}\xi$
are needed.

First, we observe that
$g(\nabla_{X}X,Z)=ag_{0}(\dot{\nabla}_{X}X,Z)$, for any $Z\in
C(T\mathbb{S}^{3})$ and for any $X\in C(T\mathbb{S}^{3})$
orthogonal to $\xi$. Then, using the properties of the Sasakian
structures $(\varphi,\xi,\eta,g)$ and
$(\varphi_{0},\xi_{0},\eta_{0},g_{0})$, we get
$$
\dot{\nabla}_{T}T=\kappa\varphi T, \ \ \dot{\nabla}_{T}\varphi
T=-\kappa T+\xi,\ \ \dot{\nabla}_{T}\xi=-\frac{1}{a}\varphi T.
$$

\noindent Further, using the Gauss equation of
$(\mathbb{S}^{3},g_{0})$ in $(\mathbb{R}^{4},\langle,\rangle)$, we
obtain
$$
\left\{\begin{array}{lll}\widetilde{\nabla}_{T}T=\dot{\nabla}_{T}T-\langle
T,T\rangle\gamma=\kappa\varphi T-\frac{1}{a}\gamma\\ \\
\widetilde{\nabla}_{T}\varphi T=-\kappa T+\xi,\ \
\widetilde{\nabla}_{T}\xi=-\frac{1}{a}\varphi T\end{array}\right..
$$

\noindent Now it follows
$$
\widetilde{\nabla}_{T}\widetilde{\nabla}_{T}T=\kappa\widetilde{\nabla}_{T}\varphi
T-\frac{1}{a}T=-\Big(\frac{1}{a}+\kappa^{2}\Big)T+\kappa\xi
$$
and
$$
\widetilde{\nabla}_{T}\widetilde{\nabla}_{T}\widetilde{\nabla}_{T}T=-\Big(\frac{2}{a}+\kappa^{2}\Big)\Big(\widetilde{\nabla}_{T}T+\frac{1}{a}\gamma\Big)+\frac{1}{a}\Big(\frac{1}{a}+\kappa^{2}\Big)\gamma.
$$

\noindent As $\widetilde{\nabla}_{T}T=\gamma''$ and
$\kappa^{2}=c-1=\frac{4}{a}-4$, we obtain the equation of $\gamma$
(thought as a curve in $\mathbb{R}^{4}$)
\begin{equation}\label{3}
a^{2}\gamma^{iv}+a(6-4a)\gamma''+\gamma=0.
\end{equation}
The general solution is
$$
\gamma(s)=\cos(As)c_{1}+\sin(As)c_{2}+\cos(Bs)c_{3}+\sin(Bs)c_{4},
$$
where $A$, $B$ are given by (\ref{2}) and $\{c_{i}\}$ are constant
vectors in $\mathbb{R}^{4}$.

It is easy to see that $\gamma$ must verify the following
relations
$$
\langle\gamma,\gamma\rangle=1,\
\langle\gamma',\gamma'\rangle=\frac{1}{a},\
\langle\gamma,\gamma'\rangle=0,\
\langle\gamma',\gamma''\rangle=0,\
\langle\gamma'',\gamma''\rangle=\frac{5-4a}{a^{2}},
$$
$$
\langle\gamma,\gamma''\rangle=-\frac{1}{a},\
\langle\gamma',\gamma'''\rangle=-\frac{5-4a}{a^{2}},\
\langle\gamma'',\gamma'''\rangle=0,\
\langle\gamma,\gamma'''\rangle=0,\
$$
$$
\langle\gamma''',\gamma'''\rangle=\frac{16a^{2}-44a+29}{a^{3}},
$$
Denoting $c_{ij}=\langle c_{i},c_{j}\rangle$, in $s=0$ we have
\begin{equation}\label{1.1}
c_{11}+2c_{13}+c_{33}=1
\end{equation}
\begin{equation}\label{1.2}
A^{2}c_{22}+2ABc_{24}+B^{2}c_{44}=\frac{1}{a}
\end{equation}
\begin{equation}\label{1.3}
Ac_{12}+Ac_{23}+Bc_{14}+Bc_{34}=0
\end{equation}
\begin{equation}\label{1.4}
A^{3}c_{12}+AB^{2}c_{23}+A^{2}Bc_{14}+B^{3}c_{34}=0
\end{equation}
\begin{equation}\label{1.5}
A^{4}c_{11}+2A^{2}B^{2}c_{13}+B^{4}c_{33}=\frac{5-4a}{a^{2}}
\end{equation}
\begin{equation}\label{1.6}
A^{2}c_{11}+(A^{2}+B^{2})c_{13}+B^{2}c_{33}=\frac{1}{a}
\end{equation}
\begin{equation}\label{1.7}
A^{4}c_{22}+(AB^{3}+A^{3}B)c_{24}+B^{4}c_{44}=\frac{5-4a}{a^{2}}
\end{equation}
\begin{equation}\label{1.8}
A^{5}c_{12}+A^{3}B^{2}c_{23}+A^{2}B^{3}c_{14}+B^{5}c_{34}=0
\end{equation}
\begin{equation}\label{1.9}
A^{3}c_{12}+A^{3}c_{23}+B^{3}c_{14}+B^{3}c_{34}=0
\end{equation}
\begin{equation}\label{1.10}
A^{6}c_{22}+2A^{3}B^{3}c_{24}+B^{6}c_{44}=\frac{16a^{2}-44a+29}{a^{3}}.
\end{equation}
Since the determinant of the system given by (\ref{1.3}),
(\ref{1.4}), (\ref{1.8}) and (\ref{1.9}) is
$-A^{2}B^{2}(A^{2}-B^{2})^{4}\neq 0$ it follows that
$$
c_{12}=c_{23}=c_{14}=c_{34}=0.
$$

The equations (\ref{1.1}), (\ref{1.5}) and (\ref{1.6}) give
$$
c_{11}=\frac{B}{A+B},\ \ c_{13}=0,\ \ c_{33}=\frac{A}{A+B},
$$
and, from (\ref{1.2}), (\ref{1.7}) and (\ref{1.10}) it results that
$$
c_{22}=\frac{B}{A+B},\ \ c_{24}=0,\ \ c_{44}=\frac{A}{A+B}.
$$
We obtained that $\{c_{i}\}$ are orthogonal vectors in
$\mathbb{R}^{4}$ with $\|c_{1}\|=\|c_{2}\|=\sqrt{\frac{B}{A+B}}$ and
$\|c_{3}\|=\|c_{4}\|=\sqrt{\frac{A}{A+B}}$. Let us consider
$c_{i}=\|c_{i}\|e_{i}$, where $\{e_{i}\}$ are mutually orthogonal
unit constant vectors in $\mathbb{R}^{4}$.

Now, if $\gamma$ is a Legendre curve, one obtains $\langle
T,\xi\rangle=0$. From the expression of $\gamma$ and since, by
definition, $\xi=-J\gamma(s)$, we get $\langle
e_{1},Je_{2}\rangle=-\langle e_{3},Je_{4}\rangle=\pm 1$ and $\langle
e_{1},Je_{3}\rangle=\langle e_{1},Je_{4}\rangle=\langle
e_{2},Je_{3}\rangle=\langle e_{2},Je_{4}\rangle=0$, from which arise
the forms of the vectors $e_{i}$. Therefore $e_{2}=\mp Je_{1}$ and
$e_{4}=\pm Je_{3}$. Finally, we obtain
$$
\begin{array}{cl}
\gamma(s)=&\sqrt{\frac{B}{A+B}}\cos(As)e_{1}-\sqrt{\frac{B}{A+B}}\sin(As)Je_{1}+\\
\\&+\sqrt{\frac{A}{A+B}}\cos(Bs)e_{3}+\sqrt{\frac{A}{A+B}}\sin(Bs)Je_{3},\end{array}
$$
and
$$
\begin{array}{cl}
\gamma_{1}(s)=&\sqrt{\frac{B}{A+B}}\cos(As)e_{1}+\sqrt{\frac{B}{A+B}}\sin(As)Je_{1}+\\
\\&+\sqrt{\frac{A}{A+B}}\cos(Bs)e_{3}-\sqrt{\frac{A}{A+B}}\sin(Bs)Je_{3},\end{array}
$$
but $\gamma$ and $\gamma_{1}$ parametrize the same curve since
$\gamma_{1}(s)=\gamma(-s)$.

\hskip11cm$\Box$

\end{pf}
\begin{remark}\textnormal{The geometric interpretation of equation
(\ref{3}) can be obtained as follows. Denote by $\textbf{j}$ the
canonical inclusion of $\gamma(\mathbb{R})$ in
$(\mathbb{S}^{3},g_{0})$. Since $\langle T,T\rangle=\frac{1}{a}$,
we obtain that
$$
\tau(\textbf{j})=a\dot{\nabla}_{T}T=a\gamma''+\gamma\ \ \hbox{and}\
\ \tau_{2}(\textbf{j})=a^{2}\gamma^{iv}+2a\gamma''+(4a-3)\gamma.
$$
Therefore $\gamma$ is a solution of (\ref{3}) if and only if
$$
\tau_{2}(\textbf{j})+4(1-a)\tau(\textbf{j})=0.
$$
We note that Riemannian immersions $\phi$ in space forms which
satisfy the equation $\tau_{2}(\phi)=\lambda\tau(\phi)$, or
equivalently, their mean curvature vectors are eigenvectors of the
rough Laplacian, were studied (for example, see \cite{Chen}).}
\end{remark}
\begin{remark}\textnormal{In \cite{Balmus} the author studied the
biharmonic curves in the Berger spheres $\mathbb{S}_{\epsilon}^{3}$.
It was proved that biharmonic curves parametrized by arc length are
helices and then their parametric equations were derived. We note
that the Berger metric is homothetic to $g$. Therefore, by changing
the parameter of the biharmonic curves of
$\mathbb{S}_{\epsilon}^{3}$ we obtain the biharmonic curves of
$(\mathbb{S}^{3},g)$, and then we may select the Legendre ones.
However, our method is completely different and the constant vectors
are precisely determined.}
\end{remark}

\section{Proper-biharmonic Hopf cylinders}
\setcounter{equation}{0}

Although proper-biharmonic curves in 3-dimensional spaces of
non-constant sectional curvature were found, various attempts to
find proper-biharmonic surfaces failed. Now, in the case of Hopf
cylinders in $(\mathbb{S}^{3},\varphi,\xi,\eta,g)$ we can state:

\begin{teo}\label{t3}
The parametric equation of the proper-biharmonic Hopf cylinder
$S_{\bar\gamma}$ in $(\mathbb{S}^{3},\varphi,\xi,\eta,g)$, thought
as a surface in $(\mathbb{R}^{4},\langle,\rangle)$, is
\begin{equation}\label{a}
\begin{array}{cl}
x=&x(u,v)=\sqrt{\frac{B}{A+B}}\cos(Au+\frac{1}{a}v)e_{1}-\sqrt{\frac{B}{A+B}}\sin(Au+\frac{1}{a}v)Je_{1}+\\ \\
&+\sqrt{\frac{A}{A+B}}\cos(Bu-\frac{1}{a}v)e_{3}+\sqrt{\frac{A}{A+B}}\sin(Bu-\frac{1}{a}v)Je_{3},
\end{array}
\end{equation}
where $\{e_{1},e_{3}\}$ is an orthonormal system of constant
vectors in the Euclidian space $(\mathbb{R}^{4},\langle,\rangle)$
with $e_{3}$ orthogonal to $Je_{1}$, and
\begin{equation}\label{b}
\left\{\begin{array}{ll}A=\sqrt{\frac{3-2a-2\sqrt{(a-1)(a-2)}}{a}}\\
\\ B=\sqrt{\frac{3-2a+2\sqrt{(a-1)(a-2)}}{a}}\end{array}\right..
\end{equation}
\end{teo}

\begin{pf} We consider the Boothby-Wang fibration $\pi:(\mathbb{S}^{3},g)\rightarrow \mathbb{C}P^{1}$, where $\mathbb{C}P^{1}$ is the complex projective space with
constant sectional holomorphic curvature $\frac{4}{a}$. Let us
denote by $\bar\nabla$, $\nabla$, $\dot{\nabla}$ and
$\widetilde{\nabla}$ the Levi-Civita connections on
$\mathbb{C}P^{1}$, $(\mathbb{S}^{3},g)$, $(\mathbb{S}^{3},g_{0})$
and $(\mathbb{R}^{4},\langle,\rangle)$, respectively.

Let $S_{\bar\gamma}$ be a proper-biharmonic Hopf cylinder in
$(\mathbb{S}^{3},g)$, where $\bar\gamma:I\rightarrow
\mathbb{C}P^{1}$ is the base curve parametrized by arc length.
Using Theorem \ref{t1.2}, $\bar\gamma$ has constant curvature
$\bar\kappa=\pm\sqrt{c-1}$.

We denote $f_{1}=\bar{\gamma}'$ and $f_{2}=\bar{J}f_{1}$. The
horizontal lift $f_{1}^{H}$ of $f_{1}$ and $\xi$ form a global
orthonormal frame field on $S_{\bar\gamma}$, and $f_{2}^{H}=\varphi
f_{1}^{H}$ is normal to $S_{\bar\gamma}$.

In order to find the explicit parametrization of $S_{\bar\gamma}$
as a surface in $\mathbb{R}^{4}$, the expressions of
$\widetilde{\nabla}_{f_{1}^{H}}f_{1}^{H}$,
$\widetilde{\nabla}_{f_{1}^{H}}\xi$,
$\widetilde{\nabla}_{f_{1}^{H}}f_{2}^{H}$ and
$\widetilde{\nabla}_{\xi}f_{2}^{H}$ are needed.

First, we have that
$$
\nabla_{f_{1}^{H}}f_{1}^{H}=(\bar{\nabla}_{f_{1}}f_{1})^{H}-g(f_{1}^{H},\varphi
f_{1}^{H})=\bar\kappa f_{2}^{H}.
$$

Using again $g(\nabla_{X}X,Z)=ag_{0}(\dot{\nabla}_{X}X,Z)$, for any
$Z\in C(T\mathbb{S}^{3})$ and for any $X\in C(T\mathbb{S}^{3})$
orthogonal to $\xi$, and the properties of the Sasakian structures
$(\varphi,\xi,\eta,g)$ and $(\varphi_{0},\xi_{0},\eta_{0},g_{0})$,
we get
$$
\dot{\nabla}_{f_{1}^{H}}f_{1}^{H}=\bar\kappa f_{2}^{H},\
\dot{\nabla}_{f_{1}^{H}}\xi=-\frac{1}{a}f_{2}^{H},\
\dot{\nabla}_{f_{1}^{H}}f_{2}^{H}=-\bar\kappa f_{1}^{H}+\xi,\
\dot{\nabla}_{\xi}f_{2}^{H}=\frac{1}{a}f_{1}^{H}.
$$
Then
$$
\widetilde{\nabla}_{f_{1}^{H}}f_{1}^{H}=\bar\kappa
f_{2}^{H}-\frac{1}{a}x,\
\widetilde{\nabla}_{f_{1}^{H}}\xi=-\frac{1}{a}f_{2}^{H},\
\widetilde{\nabla}_{f_{1}^{H}}f_{2}^{H}=-\bar\kappa
f_{1}^{H}+\xi,\
\widetilde{\nabla}_{\xi}f_{2}^{H}=\frac{1}{a}f_{1}^{H}.
$$

We recall that $[\xi,f_{1}^{H}]=0$, therefore we can choose a local
chart $x=x(u,v)$ such that $f_{1}^{H}=x_{u}$ and $\xi=x_{v}$.

After a straightforward computation, assuming that
$\bar\kappa=\sqrt{c-1}$, we obtain
\begin{equation}\label{3.01}
\left\{\begin{array}{ll}a^{2}x_{uuuu}+a(6-4a)x_{uu}+x=0\\ \\
ax_{uuv}-\sqrt{c-1}x_{u}+x_{v}=0\end{array}\right..
\end{equation}
Thus
\begin{equation}\label{3.02}
\begin{array}{cl}
x=&x(u,v)=\cos(Au+\frac{1}{a}v)c_{1}+\sin(Au+\frac{1}{a}v)c_{2}+\\
\\&+\cos(Bu-\frac{1}{a}v)c_{3}+\sin(Bu-\frac{1}{a}v)c_{4},\end{array}
\end{equation}
where $A$, $B$ are given by (\ref{b}) and $\{c_{i}\}$ are constant
vectors in $\mathbb{R}^{4}$. Since
$$
\langle x,x\rangle=1,\ \langle x,x_{u}\rangle=0,\ \langle
x,x_{v}\rangle=0,\ \langle x_{u},x_{v}\rangle=0,
$$
$$
\langle x_{u},x_{u}\rangle=\frac{1}{a},\ \langle
x_{v},x_{v}\rangle=\frac{1}{a^{2}},\ \langle
x_{v},x_{uv}\rangle=0,\ \langle x_{u},x_{uv}\rangle=0 ,
$$
one obtains, in $(u,v)=(0,0)$,
\begin{equation}\label{2.1}
c_{11}+2c_{13}+c_{33}=1
\end{equation}
\begin{equation}\label{2.2}
Ac_{12}+Bc_{14}+Ac_{23}+Bc_{34}=0
\end{equation}
\begin{equation}\label{2.3}
-c_{12}+c_{14}-c_{23}+c_{34}=0
\end{equation}
\begin{equation}\label{2.4}
-Ac_{22}+(A-B)c_{24}+Bc_{44}=0
\end{equation}
\begin{equation}\label{2.5}
A^{2}c_{22}+2ABc_{24}+B^{2}c_{44}=\frac{1}{a}
\end{equation}
\begin{equation}\label{2.6}
c_{22}-2c_{24}+c_{44}=1
\end{equation}
\begin{equation}\label{2.7}
-Ac_{12}+Ac_{14}+Bc_{23}-Bc_{34}=0
\end{equation}
\begin{equation}\label{2.8}
-A^{2}c_{12}-ABc_{14}+ABc_{23}+B^{2}c_{34}=0
\end{equation}
where $c_{ij}=\langle c_{i},c_{j}\rangle$.

From $\langle x_{u},x_{u}\rangle=\frac{1}{a}$, $\langle
x_{u},x_{v}\rangle=0$, one obtains, in $(u,v)=(0,\frac{a\pi}{2})$,
\begin{equation}\label{2.9}
A^{2}c_{11}-2ABc_{13}+B^{2}c_{33}=\frac{1}{a}
\end{equation}
\begin{equation}\label{2.10}
-Ac_{11}+(B-A)c_{13}+Bc_{33}=0.
\end{equation}
From (\ref{2.2}), (\ref{2.3}), (\ref{2.7}) and (\ref{2.8}) we have
$c_{12}=c_{14}=c_{23}=c_{34}=0$. From (\ref{2.1}), (\ref{2.9}) and
(\ref{2.10}) it follows that $c_{11}=\frac{B}{A+B}$, $c_{13}=0$,
$c_{33}=\frac{A}{A+B}$. Finally, from (\ref{2.4}), (\ref{2.5}) and
(\ref{2.6}) one obtains $c_{22}=\frac{B}{A+B}$, $c_{24}=0$,
$c_{44}=\frac{A}{A+B}$. Thus $\{c_{i}\}$ are orthogonal vectors in
$E^{4}$ with $\|c_{1}\|=\|c_{2}\|=\sqrt{\frac{B}{A+B}}$ and
$\|c_{3}\|=\|c_{4}\|=\sqrt{\frac{A}{A+B}}$. Hence
$c_{1}=\sqrt{\frac{B}{A+B}}e_{1}$,
$c_{2}=\sqrt{\frac{B}{A+B}}e_{2}$,
$c_{3}=\sqrt{\frac{A}{A+B}}e_{3}$,
$c_{4}=\sqrt{\frac{A}{A+B}}e_{4}$, where $\{e_{i}\}$ are mutually
orthogonal unit constant vectors in $\mathbb{R}^{4}$. Let us denote
by $e_{i}^{k}$, $i,k=\overline{1,4}$, the components of $e_{i}$.
Then the components of $x=x(u,v)$ are
$$
\begin{array}{cl}
x^{j}=&x^{j}(u,v)=e_{1}^{j}\cos(Au+\frac{1}{a}v)+e_{2}^{j}\sin(Au+\frac{1}{a}v)+\\
\\&+e_{3}^{j}\cos(Bu-\frac{1}{a}v)+e_{3}^{j}\cos(Bu-\frac{1}{a}v).\end{array}
$$
Since $\xi=x_{v}(u,v)$, $\xi=\frac{1}{a}\xi_{0}$, it follows that
$e_{2}=-Je_{1}$ and $e_{4}=Je_{3}$.

\noindent When $\bar\kappa=-\sqrt{c-1}$, following the same steps,
we obtain
$$
\begin{array}{cl}
x_{1}(u,v)=&\sqrt{\frac{B}{A+B}}\cos(Au-\frac{1}{a}v)e_{1}+\sqrt{\frac{B}{A+B}}\sin(Au-\frac{1}{a}v)Je_{1}+\\
\\&+\sqrt{\frac{A}{A+B}}\cos(Bu+\frac{1}{a}v)e_{3}-\sqrt{\frac{A}{A+B}}\sin(Bu+\frac{1}{a}v)Je_{3},\end{array}
$$
and we see that $x(u,v)=x_{1}(-u,v)$.

\hskip11cm$\Box$

\end{pf}
\begin{remark}
\textnormal{We fix now $\{e_{1},e_{2},e_{3},e_{4}\}$ an
orthonormal basis of $\mathbb{R}^{4}$, with $e_{2}=-Je_{1}$ and
$e_{4}=Je_{3}$, and consider $S_{\bar\gamma}$ the corresponding
Hopf cylinder. We observe that the geodesics
$$
\begin{array}{cl}
u\rightarrow
x(u,v_{0})=&\sqrt{\frac{B}{A+B}}\cos(Au+\frac{1}{a}v_{0})e_{1}+\sqrt{\frac{B}{A+B}}\sin(Au+\frac{1}{a}v_{0})e_{2}+\\
\\&+\sqrt{\frac{A}{A+B}}\cos(Bu-\frac{1}{a}v_{0})e_{3}+\sqrt{\frac{A}{A+B}}\sin(Bu-\frac{1}{a}v_{0})e_{4},\end{array}
$$
can be written as
$$
\begin{array}{cl}
\gamma(s)=&\sqrt{\frac{B}{A+B}}\cos(As)f_{1}+\sqrt{\frac{B}{A+B}}\sin(As)f_{2}+\\
\\&+\sqrt{\frac{A}{A+B}}\cos(Bs)f_{3}+\sqrt{\frac{A}{A+B}}\sin(Bs)f_{4},\end{array}
$$
where
$$
\left\{\begin{array}{lllll}f_{1}=\cos(\frac{1}{a}v_{0})e_{1}+\sin(\frac{1}{a}v_{0})e_{2}\\
\\
f_{2}=-\sin(\frac{1}{a}v_{0})e_{1}+\cos(\frac{1}{a}v_{0})e_{2}\\
\\ f_{3}=\cos(\frac{1}{a}v_{0})e_{3}-\sin(\frac{1}{a}v_{0})e_{4}\\
\\
f_{4}=\sin(\frac{1}{a}v_{0})e_{3}+\cos(\frac{1}{a}v_{0})e_{4}\end{array}\right..
$$
Note that $\{f_{1},f_{2},f_{3},f_{4}\}$ is an orthonormal basis of
$\mathbb{R}^{4}$ with $f_{2}=-Jf_{1}$ and $f_{4}=Jf_{3}$. Therefore,
the only proper-biharmonic Legendre curves of
$(\mathbb{S}^{3},\varphi,\xi,\eta,g)$ are the geodesics of the
biharmonic Hopf cylinders $S_{\bar\gamma}$ orthogonal to $\xi$.}
\end{remark}
\begin{remark}
\textnormal{In fact, $x$ is an isometric embedding $\mathcal{T}$
into $(\mathbb{S}^{3},g)$, where $\mathcal{T}$ is the flat torus
$\mathcal{T}=\mathbb{R}^{2}/\Lambda$, $\Lambda$ being the lattice
generated by the vectors
$w_{1}=\Big(\frac{2\pi}{A+B},\frac{2\pi}{A(A+B)}\Big)$ and
$w_{2}=\Big(\frac{2\pi}{A+B},-\frac{2\pi}{B(A+B)}\Big)$.}
\end{remark}
\begin{remark}\textnormal{We can obtain an alternative proof of
Theorem \ref{t3} by using Theorem \ref{t2} and the flow of the
Reeb vector field $\xi$.}
\end{remark}

A non-compact 3-dimensional Sasakian space form, locally isometric
to $(\mathbb{S}^{3},g)$, is provided by the Cartan-Vranceanu spaces.

\noindent In \cite{COP2}, the proper-biharmonic Legendre curves of
the Cartan-Vranceanu 3-dimen-sional spaces $(M^{3},ds^{2}_{l,m})$
were obtained.

\noindent We recall that $M=\mathbb{R}^{3}$ if $m\geqslant 0$ and
$M=\{(x,y,z)\in\mathbb{R}^{3}: x^{2}+y^{2}<-\frac{1}{m}\}$
otherwise, while
$$
ds^{2}_{l,m}=\frac{dx^{2}+dy^{2}}{(1+m(x^{2}+y^{2}))^{2}}+\Big(dz+\frac{l}{2}\frac{ydx-xdy}{(1+m(x^{2}+y^{2}))}\Big)^{2},
$$
$l,m\in\mathbb{R}$. When $l\neq 0$, the 1-form
$$
\eta=dz+\frac{l}{2}\frac{ydx-xdy}{(1+m(x^{2}+y^{2}))}
$$
is a contact form.

\noindent \begin{teo}[\cite{COP2}] The parametric equations of
proper-biharmonic Legendre\newline curves of $(M^{3},ds^{2}_{l,m})$,
with $4m-l^{2}>0$ and $l\neq 0$, are
\begin{equation}\label{ecop}
\left\{\begin{array}{l}x(s)=\alpha\sin(\beta s+c_1)\\ y(s)=-\alpha\cos(\beta s+c_1)\\
z(s)=\frac{l}{2}s+c_2
\end{array}\right.,
\end{equation}
where $\beta=\frac{1+m\alpha^{2}}{\alpha}$,
$\alpha^{2}=\frac{6m-l^{2}+\sqrt{32m^{2}-12ml^{2}+l^{4}}}{2m^{2}}$
and $c_1$, $c_2$ are real constants.
\end{teo}
Moreover, if $l=2$ and $m>1$ then $(\mathbb{R}^{3},ds^{2}_{2,m})$ is
a Sasakian space form with constant $\varphi$-sectional curvature
$c=4m-3>1$ (see \cite{CIL}), and $(\mathbb{R}^{3},ds^{2}_{2,m})$ is
locally isometric with $(\mathbb{S}^{3},g)$.

Therefore, the cylinders along the $Oz$-axis determined by
(\ref{ecop}), for $l=2$ and $m>1$, are the biharmonic Hopf
cylinders.

\par

\vskip1cm

\begin{tabular}{lll}
Dorel Fetcu& &Cezar Oniciuc\\
Department of Mathematics& &Faculty of Mathematics \\
Technical University of Ia\c si& &"Al. I. Cuza" University of Ia\c
si\\Bd. Carol I, no. 11& &Bd. Carol I, no. 11 \\700506 Iasi& &700506 Iasi\\
Rom\^ ania.& &
Rom\^ ania. \\
e-mail: dfetcu@math.tuiasi.ro& &e-mail: oniciucc@uaic.ro
\end{tabular}

\end{document}